\documentclass{article}
 
\usepackage{amsmath,amssymb,amsfonts,pstricks,pst-tree}

\def\Q{\mathbb{Q}}

\def\H{{\cal{H}}}
\def\Hd{{\cal{H}}^{\mathrm{dif}}}
\def\Hdnc{{\cal{H}}^{\mathrm{dif,nc}}}
\def\Ha{{\cal{H}}^{\alpha}}
\def\Hanc{{\cal{H}}^{\alpha,nc}}

\def\D{\Delta}
\def\d{\delta}
\def\Dd{\Delta^{\mathrm{dif}}}

\def\Da{\Delta^{\alpha}}

\def\DH{\Delta^{\H}}

\def\Id{\mathrm{Id}}
\def\Span{\mathrm{Span}}

\def\T{\overline{T}}

\def\calT{{\cal{T}}}

\def\over{\slash}
\renewcommand{\i}[1]{{}_{\scriptscriptstyle(#1)}}
\renewcommand{\j}[1]{{}_{\scriptscriptstyle[#1]}}

\newtheorem{theorem}{Theorem}[section]
\newtheorem{proposition}[theorem]{Proposition}
\newtheorem{corollary}[theorem]{Corollary}
\newtheorem{lemma}[theorem]{Lemma}
\newtheorem{defin}[theorem]{Definition}

\newtheorem{rem}[theorem]{Remark}
\newenvironment{remark}{\begin{rem} \em}{\end{rem}}
\newtheorem{exa}[theorem]{Examples}
\newenvironment{examples}{\begin{exa} \em}{\end{exa}}
\newenvironment{proof}{\noindent{\em Proof.\/}}{\hfill$\Box$\par\vspace{.2cm}}
\newenvironment{proof of}[2]{\bigskip \noindent{\em Proof of #1~\ref{#2}.\/}}
        {\hfill$\square$\par\vspace{.2cm}}
\numberwithin{equation}{section}
\numberwithin{figure}{section}

\newcommand{\tbc}[2]{\TBC*~[tnpos=#1]{#2}}
\newpsobject{TBC}{TC}{linecolor=black,radius=2pt}
\psset{levelsep=.5cm, treesep=.5cm,treemode=U}

\newcommand{\treeO}{
\setlength{\unitlength}{3pt}
\psset{unit=3pt}
\psset{runit=2pt}
\psset{linewidth=0.2}
\begin{pspicture}(0,0)(2,2)
\psline(1,-1)(1,1.5)
\end{pspicture}}

\newcommand{\treeA}{
\setlength{\unitlength}{3pt}
\psset{unit=3pt}
\psset{runit=2pt}
\psset{linewidth=0.2}
\begin{pspicture}(0,0)(2.5,2)
\psline(1,-1)(1,.5)
\psline(1,.5)(0,1.5)
\psline(1,.5)(2,1.5)
\end{pspicture}}

\newcommand\treeAB{
\setlength{\unitlength}{3pt}
\psset{unit=3pt}
\psset{runit=2pt}
\psset{linewidth=0.2}
\begin{pspicture}(0,0)(5,3)
\psline(3,-1)(3,.5)
\psline(3,.5)(1,2.5)
\psline(3,.5)(4,1.5)
\psline(2,1.5)(3,2.5)
\end{pspicture}}

\newcommand\treeBA{
\setlength{\unitlength}{3pt}
\psset{unit=3pt}
\psset{runit=2pt}
\psset{linewidth=0.2}
\begin{pspicture}(0,0)(5,3)
\psline(2,-1)(2,.5)
\psline(2,.5)(4,2.5)
\psline(2,.5)(1,1.5)
\psline(3,1.5)(2,2.5)
\end{pspicture}}

\newcommand\treeABC{
\setlength{\unitlength}{3pt}
\psset{unit=3pt}
\psset{runit=2pt}
\psset{linewidth=0.2}
\begin{pspicture}(0,0)(5,4.5)
\psline(3,-1)(3,.5)
\psline(3,.5)(0,3.5)
\psline(3,.5)(4,1.5)
\psline(2,1.5)(3,2.5)
\psline(1,2.5)(2,3.5)
\end{pspicture}}

\newcommand\treeBAC{
\setlength{\unitlength}{3pt}
\psset{unit=3pt}
\psset{runit=2pt}
\psset{linewidth=0.2}
\begin{pspicture}(0,0)(5,4.5)
\psline(3,-1)(3,.5)
\psline(3,.5)(1,2.5)
\psline(3,.5)(4,1.5)
\psline(2,1.5)(4,3.5)
\psline(3,2.5)(2,3.5)
\end{pspicture}}

\newcommand\treeACA{
\setlength{\unitlength}{3pt}
\psset{unit=3pt}
\psset{runit=2pt}
\psset{linewidth=0.2}
\begin{pspicture}(0,0)(6,3.5)
\psline(3,-1)(3,.5)
\psline(3,.5)(0.5,3)
\psline(3,.5)(5.5,3)
\psline(1.5,2)(2.5,3)
\psline(4.5,2)(3.5,3)
\end{pspicture}}

\newcommand\treeCAB{
\setlength{\unitlength}{3pt}
\psset{unit=3pt}
\psset{runit=2pt}
\psset{linewidth=0.2}
\begin{pspicture}(0,0)(5,4.5)
\psline(2,-1)(2,.5)
\psline(2,.5)(4,2.5)
\psline(2,.5)(1,1.5)
\psline(3,1.5)(1,3.5)
\psline(2,2.5)(3,3.5)
\end{pspicture}}

\newcommand\treeCBA{
\setlength{\unitlength}{3pt}
\psset{unit=3pt}
\psset{runit=2pt}
\psset{linewidth=0.2}
\begin{pspicture}(0,0)(5,4.5)
\psline(2,-1)(2,.5)
\psline(2,.5)(5,3.5)
\psline(2,.5)(1,1.5)
\psline(3,1.5)(2,2.5)
\psline(4,2.5)(3,3.5)
\end{pspicture}}

\newcommand\treeABDA{
\setlength{\unitlength}{3pt}
\psset{unit=3pt}
\psset{runit=2pt}
\psset{linewidth=0.2}
\begin{pspicture}(0,0)(7,4.5)
\psline(4,-1)(4,.5)
\psline(4,.5)(.5,4)
\psline(4,.5)(6.5,3)
\psline(2.5,2)(3.5,3)
\psline(1.5,3)(2.5,4)
\psline(5.5,2)(4.5,3)
\end{pspicture}}

\newcommand\treeADAB{
\setlength{\unitlength}{3pt}
\psset{unit=3pt}
\psset{runit=2pt}
\psset{linewidth=0.2}
\begin{pspicture}(0,0)(6,4.5)
\psline(3,-1)(3,.5)
\psline(3,.5)(.5,3)
\psline(3,.5)(5.5,3)
\psline(1.5,2)(2.5,3)
\psline(4.5,2)(2.5,4)
\psline(3.5,3)(4.5,4)
\end{pspicture}}

\newcommand\treeDABC{
\setlength{\unitlength}{3pt}
\psset{unit=3pt}
\psset{runit=2pt}
\psset{linewidth=0.2}
\begin{pspicture}(0,0)(5.5,5)
\psline(3,-1)(3,.5)
\psline(3,.5)(2,1.5)
\psline(3,.5)(5,2.5)
\psline(4,1.5)(1,4.5)
\psline(3,2.5)(4,3.5)
\psline(2,3.5)(3,4.5)
\end{pspicture}}

\newcommand\treeDACA{
\setlength{\unitlength}{3pt}
\psset{unit=3pt}
\psset{runit=2pt}
\psset{linewidth=0.2}
\begin{pspicture}(0,0)(6.5,5)
\psline(3,-1)(3,.5)
\psline(3,.5)(2,1.5)
\psline(3,.5)(6.5,4)
\psline(4,1.5)(1.5,4)
\psline(2.5,3)(3.5,4)
\psline(5.5,3)(4.5,4)
\end{pspicture}}

\newcommand\treeDBAC{
\setlength{\unitlength}{3pt}
\psset{unit=3pt}
\psset{runit=2pt}
\psset{linewidth=0.2}
\begin{pspicture}(0,0)(5.5,5)
\psline(2,-1)(2,.5)
\psline(2,.5)(1,1.5)
\psline(2,.5)(4,2.5)
\psline(3,1.5)(1,3.5)
\psline(2,2.5)(4,4.5)
\psline(3,3.5)(2,4.5)
\end{pspicture}}

\newcommand\treeDCAB{
\setlength{\unitlength}{3pt}
\psset{unit=3pt}
\psset{runit=2pt}
\psset{linewidth=0.2}
\begin{pspicture}(0,0)(5.5,5)
\psline(2,-1)(2,.5)
\psline(2,.5)(1,1.5)
\psline(2,.5)(5,3.5)
\psline(3,1.5)(2,2.5)
\psline(4,2.5)(2,4.5)
\psline(3,3.5)(4,4.5)
\end{pspicture}}

\newcommand\treeEABDA{
\setlength{\unitlength}{3pt}
\psset{unit=3pt}
\psset{runit=2pt}
\psset{linewidth=0.2}
\begin{pspicture}(0,0)(7,5.5)
\psline(3,-1)(3,.5)
\psline(3,.5)(2,1.5)
\psline(3,.5)(6.5,4)
\psline(4,1.5)(.5,5)
\psline(2.5,3)(3.5,4)
\psline(1.5,4)(2.5,5)
\psline(5.5,3)(4.5,4)
\end{pspicture}}

\newcommand\treeEADAB{
\setlength{\unitlength}{3pt}
\psset{unit=3pt}
\psset{runit=2pt}
\psset{linewidth=0.2}
\begin{pspicture}(0,0)(6,5.5)
\psline(2,-1)(2,.5)
\psline(2,.5)(1,1.5)
\psline(2,.5)(5.5,4)
\psline(3,1.5)(.5,4)
\psline(1.5,3)(2.5,4)
\psline(4.5,3)(2.5,5)
\psline(3.5,4)(4.5,5)
\end{pspicture}}

\newcommand\treeEDABC{
\setlength{\unitlength}{3pt}
\psset{unit=3pt}
\psset{runit=2pt}
\psset{linewidth=0.2}
\begin{pspicture}(0,0)(5.5,6)
\psline(2,-1)(2,.5)
\psline(2,.5)(1,1.5)
\psline(2,.5)(5,3.5)
\psline(3,1.5)(2,2.5)
\psline(4,2.5)(1,5.5)
\psline(3,3.5)(4,4.5)
\psline(2,4.5)(3,5.5)
\end{pspicture}}

\newcommand{\lrgraft}[2]
{\setlength{\unitlength}{3pt}
\psset{unit=3pt}
\psset{runit=2pt}
\psset{linewidth=0.2}
\begin{pspicture}(0,0)(6,3.5)
\psline(3,-1)(3,.5)
\psline(3,.5)(2,1.5)
\psline(3,.5)(4,1.5)
\put(0.5,3.5){$#1$}
\put(6.5,3.5){$#2$}
\end{pspicture}}

\newcommand{\lgraft}[2]
{\setlength{\unitlength}{4pt}
\psset{unit=5pt}
\psset{runit=4pt}
\psset{linewidth=.1}
\begin{pspicture}(0,0)(4,4)
\psline(2.1,.9)(1.4,1.6)
\put(2.6,0){$#1$}
\put(0.5,2.6){$#2$}
\end{pspicture}}

\newcommand{\rgraft}[2]
{\setlength{\unitlength}{4pt}
\psset{unit=5pt}
\psset{runit=4pt}
\psset{linewidth=.1}
\begin{pspicture}(0,0)(4,4)
\psline(1.5,.9)(2.2,1.6)
\put(.4,0){$#1$}
\put(3,2.6){$#2$}
\end{pspicture}}

\newcommand{\combRgraft}[3]
{\setlength{\unitlength}{5pt}
\psset{unit=6pt}
\psset{runit=5pt}
\psset{linewidth=.1}
\begin{pspicture}(-1,.5)(9,7)
\psline(3,-.5)(3,1)
\psline(3,1)(8,6)
\psline(3,1)(2,2)
\psline(4,2)(3,3)
\psline(7,5)(6,6)
\put(4.4,4.4){.}
\put(5.2,5.2){.}
\put(6,6){.}
\put(0,3){$#1$}
\put(1.5,4.7){$#2$}
\put(5.2,8){$#3$}
\end{pspicture}}

\begin{document}
 
\title{Combinatorial Hopf algebras from renormalization}
 
\author{Christian Brouder \\ 
Institut de Min\'eralogie et de Physique des Milieux Condens\'es, \\ 
CNRS UMR7590, Universit\'es Paris 6 et 7, IPGP, \\
140 rue de Lourmel, 75015 Paris, France. \cr
{\small \tt christian.brouder@impmc.jussieu.fr}
\and
Alessandra Frabetti \cr
Universit\'e de Lyon, Universit\'e Lyon 1, CNRS,  
UMR 5208 Institut Camille Jordan, \\ 
43, blvd du 11 novembre 1918, 69622 Villeurbanne, France. \cr
{\small \tt frabetti@math.univ-lyon1.fr}
\and
Fr\'ed\'eric Menous \\ 
D\'epartement de Math\'ematiques, B\^at. 425, \\ 
Universit\'e Paris-Sud, 91405 Orsay, France. \cr
{\small \tt Frederic.Menous@math.u-psud.fr} 
}

\date{September 18, 2009}
 
\maketitle

%%%%%%%%%%%%%%%%%%%%%%%%%%%%%%%%%%%%%%%%%%%%%%%%%%%%%%%%%%%%%%%%%%%%

\begin{abstract}
In this paper we describe the right-sided combinatorial Hopf structure of
three Hopf algebras appearing in the context of renormalization in
quantum field theory: the non-commutative version of the Fa\`a di
Bruno Hopf algebra, the non-commutative version of the charge 
renormalization Hopf algebra on planar binary trees for quantum 
electrodynamics, and the non-commutative version of the Pinter 
renormalization Hopf algebra on any bosonic field. 

We also describe two general ways to define the associative 
product in such Hopf algebras, the first one by recursion, and the
second one by grafting and shuffling some decorated rooted trees. 
\end{abstract}

%%%%%%%%%%%%%%%%%%%%%%%%%%%%%%%%%%%%%%%%%%%%%%%%%%%%%%%%%%%%%%%%%%%%

%\tableofcontents

%%%%%%%%%%%%%%%%%%%%%%%%%%%%%%%%%%%%%%%%%%%%%%%%%%%%%%%%%%%%%%%%%%%%

\section{Introduction}

In the paper \cite{LodayRonco2008}, J.-L.~Loday and M.~Ronco gave the 
definition of a right-sided combinatorial Hopf algebra (CHA) which 
includes the examples of Hopf algebras describing the 
renormalization in quantum field theory. 
The toy example, which inspired the work in \cite{LodayRonco2008}, 
is Kreimer's Hopf algebra on rooted trees introduced in \cite{Kreimer}. 
In this paper, we describe the CHA structure of three other 
such examples: the non-commutative version of the Fa\`a di Bruno Hopf algebra 
given in \cite{BFK}, the non-commutative version of the charge 
renormalization Hopf algebra on planar binary trees given in 
\cite{BFqedren} for quantum electrodynamics, and the non-commutative 
version of the Pinter renormalization Hopf algebra on any quantum
fields, as given in \cite{BrouderSchmitt}. 

Since Loday and Ronco work in the context of cofree-coassociative CHAs,
we describe our examples, which are all free-associative CHAs, by
considering their linear duals. 
A right-sided cofree-coassociative CHA is completely determined by 
the brace structure on the set its primitive elements. Therefore, 
to describe our examples, it will be sufficient to give their
brace structure. 

In the first section we fix the notation on right-sided 
cofree-coassociative CHAs and brace algebras from \cite{LodayRonco2008}. 
We also give a recursive definition of the product and show how to 
construct it starting from the one defined on decorated rooted trees, 
that is, the product defined on the right-sided CHA induced by the 
free brace algebra over the space of decorations. 
Then we fix the notation for the dual Hopf algebras, which are
used in the sequel. 
Each example is then treated in a separate section. 

% \todo{THIS IS SUSPENDED BECAUSE I STILL DIDN'T FIND THE PROOF. 

% Beside the explicit formulas for the brace products, in the third
% section we prove that the main example of a right-sided
% cofree-coassociative CHA given in \cite{LodayRonco2008}, namely the
% dendriform algebra of planar binary trees introduced in
% \cite{LodayRonco1998}, is in fact the dual of the non-commutative
% charge renormalization Hopf algebra $\Ha$ given in \cite{BFqedtree}. 
% This adds a new proof to the two already existing ones,
% cf. \cite{Foissy,Holtkamp}. } 

%%%%%%%%%%%%%%%%%%%%%%%%%%%%%%%%%%%%%%%%%%%%%%%%%%%%%%%%%%%%%%%%%%%%

\section{Combinatorial Hopf algebras from 
dual Hopf algebras}
\label{section1}

\subsection{Right-sided combinatorial Hopf algebras}
\label{rs-CHA}

In this section we recall the definition and main properties of a 
right-sided combinatorial Hopf algebra (r-s CHA), as given by J.-L.~Loday and 
M.~Ronco in \cite{LodayRonco2008}, in the cofree-coassociative case. 

A right-sided cofree-coassociative CHA is, up to an isomorphism, 
a tensor coalgebra $T(R)=\oplus_{n\geq 0} R^{\otimes n}$ with the 
deconcatenation coproduct, endowed with an associative product $\star$
which makes it into a bialgebra and which satisfies the right-sided condition: 
for any $p\geq 0$, the subspace 
$T^{\geq p}(R) = \oplus_{n\geq p} R^{\otimes n}$ 
is a right ideal of $T(R)$, that is
\begin{align*}
T^{\geq p}(R) \star T(R) \quad\subset\quad T^{\geq p}(R). 
\end{align*}
This condition is equivalent to demanding that the product $\star$ induces 
a right action of $T(R)$ on the quotient space 
$T_{< p}(R)= T(R)/T^{\geq p}(R)$, that is 
\begin{align}
\label{right-sided-cofree}
T_{< p}(R) \star T(R) \quad\subset\quad T_{< p}(R). 
\end{align}

Loday and Ronco proved that a cofree-coassociative CHA is right-sided 
if and only if the set of its primitive elements $R$ is a brace
algebra, that is, it is endowed with a brace product $\{\ ;\, ...,\ \}$
acting on $R\otimes T(R)$ with value in $R$, satisfying the brace relation 
\begin{align}
\label{brace relation}
\big\{\{x;y_1,\ldots,y_n\};z_1,\ldots,z_m\big\} 
&= \sum\ \big\{x;z_1,\ldots,\{y_1;z_{k_1},\ldots,\},  
\ldots,\{y_n;z_{k_n},\ldots \},\ldots,z_m\big\},  
\end{align}
for any $x,y_i,z_j \in R$. 
The brace product restricted to $R\otimes R^{\otimes q}$ is also denoted by 
$M_{1q}$, if one wishes to specify the number of variables on the
right, or to underline that it is a special case of a multibrace 
product $M_{pq}$ defined on any powers $R^{\otimes p}\otimes R^{\otimes q}$. 

The brace product can be found from the associative product $\star$ by
projecting the result onto the space $R$ of co-generators. 
If we denote by $\pi:T(R)\longrightarrow R$ the canonical projection, 
this means that 
\begin{align*}
M_{1q}(x;y_1\cdots y_q) &= \pi\big(x\star(y_1\cdots y_q)\big),  
\end{align*}
for any $x\in R$ and $y_1\cdots y_q\in R^{\otimes q}$, where we denote
the tensor product in $T(R)$ by the concatenation. 
In particular, for any $x\in R$ we have 
\begin{align*}
M_{10}(x;1)\equiv M_{10}(x) &= \pi(x\star 1) = x, 
\end{align*}
therefore $M_{10}=\Id$ on $R$. 

Conversely, the product $\star$ in $T(R)$ can be reconstructed 
from the brace product on $R$, using the deconcatenation coproduct to
recover $T(R)$ from $R$. 
The resulting formula of Loday and Ronco for the product is 
(cf.~\cite{LodayRonco2008}) 
\begin{align}
\label{right-sided product}
(x_1\cdots x_n)\star (y_1\cdots y_m) 
&= \sum_{k\geq 1}\ 
\underset{p_1,...,p_k=0,1}{\sum_{p_1+\cdots+p_k=n}}\  
\underset{q_1,...,q_k\geq 0}{\sum_{q_1+\cdots+q_k=m}}\  
[M_{p_1 q_1}\cdots M_{p_k q_k}](x_1\cdots x_n;y_1\cdots y_m) , 
\end{align}
where each map $M_{pq}$ is applied to a block of $p$ generators $x$
and $q$ generators $y$, with the following assumptions: 
\begin{align}
\label{brace-assumption}
& M_{00}=0 \qquad\mbox{and}\qquad M_{01}=M_{10}=\Id, \nonumber \\ 
& M_{0p}=M_{p0}=0 \qquad\mbox{for $p>1$}, \\ 
& M_{pq}=0 \qquad\mbox{for $p>1$ and $q\neq 0$}, \nonumber
\end{align}
and where the operators $M_{1q}$ for $q\neq 0$ of course satisfy 
the brace relation (\ref{brace relation}). 
\bigskip 

Applying formula (\ref{right-sided product}) requires some
computations. Let us take the example of $x \star y$.
For $k=1$ we have $p_1=q_1=1$ and the only term is
$M_{11}$. For $k=2$ we have 
$(p_1,p_2)=(1,0)$ or $(p_1,p_2)=(0,1)$ and
$(q_1,q_2)=(1,0)$ or $(q_1,q_2)=(0,1)$.
This gives us the four terms
$M_{11}M_{00}+M_{10}M_{01}+M_{01}M_{10}+M_{00}M_{11}$.
Since $M_{00}=0$, we have only two terms left
$M_{10}M_{01}+M_{01}M_{10}$.
For $k>2$, all the terms are zero because each term contains
a factor $M_{00}$.
In general, the terms with $k> n+m$ are zero in the product
$(x_1\cdots x_n)\star (y_1\cdots y_m)$.
Therefore, using (\ref{brace-assumption}), we get
\begin{align*}
x\star y &= M_{11}(x;y)+[M_{10}M_{01}](x;y)+[M_{01}M_{10}](x;y) \\
&= M_{11}(x;y) +M_{10}(x) M_{01}(y)+M_{01}(y)M_{10}(x)
= M_{11}(x;y) + xy+ yx. 
\end{align*}

Similarly, for $(x y) \star z$, we have two terms for $k=2$:
$M_{11}M_{10}+M_{10}M_{11}$
and three terms for $k=3$:
$M_{01}M_{10}M_{10}+M_{10}M_{01}M_{10}+M_{10}M_{10}M_{01}$, 
so that
\begin{eqnarray*}
(x y)\star z &=& x M_{11}(y;z) + M_{11}(x;z)y + zxy + xzy + xyz.
\end{eqnarray*}

For $x\star (yz)$, we have one term for $k=1$: $M_{12}$, 
two terms for $k=2$:
$M_{11}M_{01}+M_{01}M_{11}$
and three terms for $k=3$:
$M_{01}M_{01}M_{10}+M_{01}M_{10}M_{01}+M_{10}M_{01}M_{01}$, 
so that
\begin{eqnarray*}
x \star (yz) &=& M_{12}(x;yz)+yM_{11}(x;z) + M_{11}(x;y)z + yzx+yxz+xyz.
\end{eqnarray*}

In particular then, formula (\ref{right-sided product}) for $n=1$ gives 
\begin{align}
\label{braceexample}
x\star (y_1\cdots y_m) 
&= \sum_{i=0}^m \sum_{q=0}^{m-i}\  
y_1\cdots y_{i}\ M_{1q}(x;y_{i+1}\cdots y_{i+q})\ y_{i+q+1}\cdots y_m. 
\end{align}

%%%%%%%%%%%%%%%%%%%

\subsection{Recursive definition of the product.}
\label{sec:formulas}

Formula (\ref{right-sided product}) can be given in a recursive way. 
To do it, let us denote by $M:R\otimes T(R)\longrightarrow R$ 
the brace product $\{\ ;\ ,\ldots,\ \}$ induced on $R$ by the
associative product $\star$ on $T(R)$, and let us denote by 
$X,Y,Z...$ the generic words in $T(R)$. 

\begin{proposition}
Given the brace product $M$ on $R$, the $\star$ product on $T(R)$ 
can be reconstructed in the following recursive way. 
For any $X,Y\in T(R)$, with $X,Y \notin T^0(R)$:  
\begin{align}
\label{product-recursive}
X\star Y &= \sum_{X = X^1 X^2 \atop Y = Y^1 Y^2} (X^1\star Y^1) M(X^2;Y^2)\\
&= \sum_{X = X^1 X^2 \atop Y = Y^1 Y^2} M(X^1;Y^1) (X^2 \star Y^2) 
\nonumber \\
&= \sum_{X = X^1 X^2 X^3\atop Y = Y^1 Y^2 Y^3} 
(X^1 \star Y^1) M(X^2, Y^2) (X^3 \star Y^3), \nonumber
\end{align} 
where the sums run over all possible factorisations of $X$ and $Y$ 
with respect to the concatenation. 
\end{proposition}

\begin{proof}
With the assumptions (\ref{brace-assumption}) on $M$, formula 
(\ref{right-sided product}) can be expressed as
\begin{align}
\label{rs} 
X \star Y &= 
\sum_{k \geq 1} \sum_{X=X^1\cdots X^k \atop Y=Y^1\cdots Y^k} 
M(X^1;Y^1)\cdots M(X^k;Y^k), 
\end{align}
where the sum runs other all the factorisations (for the
concatenation) of $X=X^1 \cdots X^k$ and $Y=Y^1 \cdots Y^k$. 
Note that the factors can be in $T^0(R)$ but the condition $M(1;1)=0$
ensures that the sum is finite. 
It also ensures that formula (\ref{product-recursive}) is recursive. 

If we then use this formula to express the factor $X^1\star Y^1$ in line 
(\ref{product-recursive}), we obtain the desired expression for 
$X\star Y$. 
\end{proof}

In particular, if we expand the words $X$ and $Y$ of $T(R)$ in terms 
of elements $x_i$ and $y_j$ of $R$, with the notation of 
Section~\ref{rs-CHA}, formula (\ref{product-recursive}) becomes 
\begin{align}
\label{rs-product-recursive}
(x_1\cdots x_n)\star (y_1\cdots y_m) 
&= \big[(x_1\cdots x_{n-1})\star (y_1\cdots y_m)\big] x_n
+ \big[(x_1\cdots x_n)\star (y_1\cdots y_{m-1})\big] y_m \nonumber \\ 
&\qquad + \sum_{i=1}^{m-1}  \big[(x_1\cdots x_{n-1})\star (y_1\cdots y_i)\big]
M_{1 m-i}(x_n;y_{i+1}\cdots y_m).  
\end{align}

\begin{corollary}
The brace relation (\ref{brace relation}) is equivalent to 
the identity 
\begin{align}
\label{brace-recursive}
M(M(X;Y);Z)&=M(X;Y\star Z), 
\end{align}
for any $X\in R$ and $Y,Z \in T(R)$. 
\end{corollary}

\begin{proof}
If we call $X=x$, $Y=y_1\cdots y_n$ and $Z=z_1\cdots z_m$, the left
hand-side of the brace relation (\ref{brace relation}) is simply 
$M(M(X;Y);Z)$. The right-hand side, instead, becomes 
\begin{align*}
\sum_{k \geq 1} \sum_{Y=Y^1\cdots Y^k \atop Z=Z^1\cdots Z^k} 
M(X;M(Y^1;Z^1)\cdots M(Y^k;Z^k)), 
\end{align*}
where we use the assumptions (\ref{brace-assumption}) on $M$. 
From (\ref{rs}) it then follows that the right-hand side is equal to 
$M(X;Y\star Z)$. 
\end{proof}

\begin{remark}
The right-sided combinatorial Hopf algebras, as well as the shuffle 
and the quasishuffle Hopf algebras, are examples of more general 
cofree-coassociative combinatorial Hopf algebras $T(R)$ introduced 
by Loday and Ronco, where $R$ is a so-called multibrace algebra. 
See \cite{LodayRonco2008} and the references therein for details. 
In all these examples, as well as in the general case, the product
$\star$ on $T(R)$ can be defined recursively from the multibrace
product on $R$. The expressions (\ref{rs}), (\ref{product-recursive}) and 
(\ref{brace-recursive}), with $X\in T(R)$, are valid in the most 
general case as well as in all the examples. 
\end{remark}

%%%%%%%%%%%%%%%%%%%%%%%%%%%%%

\subsection{Right-sided Hopf algebras and trees}
\label{sec:trees}

Let $T(R)$ be a right-sided CHA such that the brace algebra $R$ is a 
finite dimensional or a graded vector space. 
In the first part of this section we forget the brace structure on $R$
and consider the free brace algebra on $R$ modeled by planar rooted trees. 
The free brace algebra has already been described in terms of trees by 
F.~Chapoton in \cite{Chapoton}, as a consequence of his construction of 
the Brace operad. Here we give a direct description of the free brace 
algebra within the notations we established in section \ref{sec:formulas}. 

In the second part of this section we relate the free brace structure on $R$ 
with the original one, and deduce a way to describe the associative product
$\star$ on $T(R)$ starting from the one defined on trees. Again, this is 
a particular case of the equivalence between the categories of brace and 
dendriform algebras given by Chapoton in \cite{Chapoton}. However it seems 
to us useful for computations to write an explicit statement. 
\bigskip

Let $\calT$ denote the set of plane rooted trees and for any non-empty
tree $t$ denote by $V(t)$ the set of its vertices and by $|t|$ the 
cardinality of $V(t)$. 
If $d$ is a map from $V(t)$ to $R^{\otimes |t|}$, then the pair
$(t,d)$ is a plane tree with decorations in $R$. For example, the 
following tree with 4 vertices is decorated by the given map $d$:
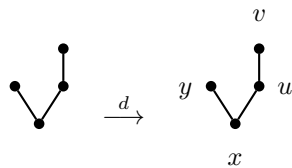
$$
\pstree{\TBC*}{\TBC* \pstree{\TBC*}{\TBC*}} 
\quad\stackrel{d}{\longrightarrow}\quad
\pstree{\tbc{b}{x}}{\tbc{l}{y} \pstree{\tbc{r}{u}}{\tbc{a}{v}}} 
$$

For any $t\in\calT$, we denote by $R^t$ the linear span of all the 
decorated trees of shape $t$, modulo the obvious linear relations 
coming from $R$, such as:
$$
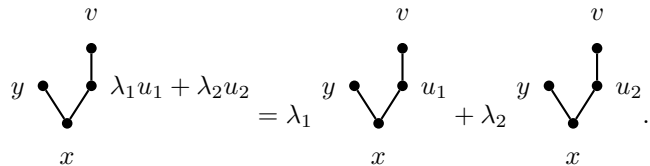
\pstree{\tbc{b}{x}}{\tbc{l}{y} \pstree{\tbc{r}{\lambda_1 u_1 +\lambda_2
    u_2}}{\tbc{a}{v}}}= \lambda_1 \ \pstree{\tbc{b}{x}}{\tbc{l}{y}
  \pstree{\tbc{r}{u_1}}{\tbc{a}{v}}}  +  \lambda_2 \ 
\pstree{\tbc{b}{x}}{\tbc{l}{y} \pstree{\tbc{r}{u_2}}{\tbc{a}{v}}}  . 
$$
Then the vector space 
$$
R^{\calT}=\bigoplus_{t\in \calT} R^t
$$
contains the decorated trees of any shape, and the elements of the 
tensor space $T(R^{\calT})$ are called decorated forests. 
From now on, we denote by the concatenation the tensor product between
trees, and we simply denote by $t$ a decorated tree, omitting
the decoration.
\bigskip 

We now define the structure of a brace algebra on $R^{\calT}$, and 
consequently that of a right-sided CHA on $T(R^{\calT})$, using the
results of section~\ref{sec:formulas}. 
Therefore, we endow $T(R^{\calT})$ with the free (deconcatenation) coproduct 
and we suppose that the product $\star$ is defined through the brace 
product $M$ on $R^{\calT}$. 

To define the brace product, consider, for any $x\in R$, the linear map 
\begin{align*}
B_+^x: T(R^{\calT}) \longrightarrow R^{\calT}
\end{align*}
which maps a forest $f=t_1\cdots t_n$ into the tree $t=B_{+}^{x}(f)$ 
obtained by joining the ordered trees $t_1,...,t_n$ to a new root 
decorated by $x$. 
For instance
$$
B_{+}^{x} \left (
  \pstree{\tbc{b}{y}}{} \hspace{.3cm}
  \pstree{\tbc{b}{u}}{\tbc{a}{v}} \right )
= \pstree{\tbc{b}{x}}{\tbc{l}{y}
  \pstree{\tbc{r}{u}}{\tbc{a}{v}}}. 
$$
Any tree $t$ can be written as $B_{+}^{x}(f)$. 
If $t$ is made only of the root, decorated by $x$, the forest $f$ is the
formal unit of $T(R^{\calT})$, that is $1\in (R^{\calT})^0$. 

\begin{theorem} 
Let $M:R^{\calT}\otimes T(R^{\calT})\longrightarrow R^{\calT}$ be the linear
map defined on any decorated tree $t=B_{+}^{x}(f_1)$ and any decorated
forest $f_2$ by 
\begin{align}
\label{def-brace-trees}
M(B_{+}^{x}(f_1);f_2) &= B_{+}^{x}(f_1 \star f_2), 
\end{align}
where $\star$ is the product induced by $M$ on $T(R^{\calT})$ according 
to Eq.~(\ref{rs}) and with the assumptions (\ref{brace-assumption}). 

Then $M$ is a brace product on $R^{\calT}$, and consequently 
$\star$ is an associative product on $T(R^{\calT})$ which makes it into a
right-sided CHA. 
\end{theorem}

\begin{proof}
Assuming conditions (\ref{brace-assumption}) on the maps $M_{pq}$, it 
suffices to show that $M$ satisfies the brace 
identity~(\ref{brace-recursive}), that is 
$$
M(M(t_1;f_2);f_3)=M(t_1;f_2\star f_3)
$$
for any tree $t_1$ and any forests $f_2,f_3$. 
If $t_1$, $f_2$, $f_3$ have respectively $n_1$, $n_2$ and $n_3$ 
vertices, we prove the brace identity by induction on the total number 
of vertices $n_1+n_2+n_3$. 

Since $t_1$ is a tree, it has at least one vertex (the root) and therefore 
$n_1+n_2+n_3\geq 1$. 
If $n_1+n_2+n_3= 1$, only $t_1$ is a tree, namely the single root, 
and $f_2,f_3$ are scalars. Because $1\star 1=1$ and $M_{10}=M_{01}=\Id$, 
both terms in the equality produce $t_1$ and therefore coincide. 

Now suppose that the brace identity holds for $n_1+n_2+n_3\leq n$, 
for a given $n\geq 1$. 
According to the previous sections, this implies that the product 
$\star$ is associative up to the same number $n$ of vertices. 
(In fact, the associativity of $\star$ does not depend on the
definition (\ref{def-brace-trees}) of $M$ and was proved in general by
Loday and Ronco.) 
Then consider a tree and two forests with a total number
of vertices $n_1+n_2+n_3= n+1$. 
Using the definition of $M$ and the associativity of $\star$ up to $n$, 
for $t_1=B_{+}^x(f_1)$ we obtain 
\begin{align*}
M(M(B_{+}^x(f_1); f_2);f_3) &= M(B_{+}^x(f_1 \star f_2);f_3) \\ 
&= B_{+}^x((f_1 \star f_2)\star f_3)  \qquad\text{(Recursion)} \\
&= B_{+}^x(f_1\star (f_2\star f_3)) \\ 
&= M(B_{+}^x(f_1);f_2\star f_3). 
\end{align*}
\end{proof}
  
\begin{lemma} 
The brace product $M$ on $R^{\calT}$ can be expressed in terms of graftings
of trees: if $t\in R^{\calT}$ and $s_1\cdots s_n\in T(R^{\calT})$, then 
$M(t;s_1\cdots s_n)$ is the sum of the terms obtained by grafting the 
trees $s_1,...s_n$ onto the vertices of $t$, in all the possible ways
which preserve the order of the $s_i$'s from left to right. 

Similarly, the product $\star$ on $T(R^{\calT})$ is the sum of terms 
mixing all possible shuffles and all possible graftings of the
trees on the right-hand side onto the vertices of the trees on the 
left hand-side.
\end{lemma}

This result shows that the recursive definition of the brace product 
(\ref{def-brace-trees}) is equivalent to the definition given in 
\cite{Chapoton} in terms of grafting of trees. 
\bigskip 

\begin{proof}
We prove it by induction on the size of the forests, given by the total
number of vertices. For this, we need to treat at the same time $M$
and $\star$, which is given by (\ref{product-recursive}). 

The induction starts with the first non-trivial operations: 
$$
M\left(\pstree{\tbc{b}{x}}{};\pstree{\tbc{b}{y}}{}\right) = B_+^x
\left (\pstree{\tbc{b}{y}}{} \right ) = \pstree{\tbc{b}{x}}{\tbc{a}{y}} 
$$ 
is the grafting of the second $\bullet$ on the first one. 
The claim is obtained because this is the only possible grafting. 
Similarly, according to (\ref{product-recursive}), the product 
$$
\pstree{\tbc{b}{x}}{} \star \pstree{\tbc{b}{y}}{} =
M\left (\pstree{\tbc{b}{x}}{};\pstree{\tbc{b}{y}}{}\right ) +
\left(\pstree{\tbc{b}{x}}{} \star 1\right )
M\left(1;\pstree{\tbc{b}{y}}{}\right  ) 
+ \left(1\star \pstree{\tbc{b}{y}}{}\right  )
M\left(\pstree{\tbc{b}{x}}{};1\right  ) 
= \pstree{\tbc{b}{x}}{\tbc{a}{y}}+ \pstree{\tbc{b}{x}}{} \hspace{.3cm}
\pstree{\tbc{b}{y}}{} + \pstree{\tbc{b}{y}}{} \hspace{.3cm}
\pstree{\tbc{b}{x}}{} 
$$ 
is exactly made of the only possible grafting and the two allowed shuffles. 

Now suppose that $M(t_1;f_2)$ is made of all possible ordered graftings of the
trees composing the forest $f_2$ on the vertices of the tree $t_1$,
for $t_1$ and $f_2$ having a total number of vertices $n_1+n_2\leq n$,
for some $n\geq 1$. 
Similarly, suppose that $f_1\star f_2$ is made of the terms mixing all
possible shuffles of the trees of $f_1$ with those of $f_2$, and of
graftings of the trees of $f_2$ over the vertices of those contained
in $f_1$, for $f_1$ and $f_2$ having a total number of vertices 
$n_1+n_2\leq n$. 

Then, consider $t_1$ and $f_2$ with $n_1+n_2=n+1$. Applying the 
definition~(\ref{def-brace-trees}) for $t_1=B_+^x(f_1)$, we see that 
$$
M(B_+^x(f_1);f_2) = B_+^x(f_1\star f_2) 
$$ 
is the grafting of $f_1\star f_2$ onto the original root. The claim 
is then guaranteed by the inductive hypothesis on $f_1\star f_2$. 

Similarly, consider $t_1$ and $f_2$ with $n_1+n_2=n+1$. Applying 
formula~(\ref{product-recursive}) we see that 
$$
f_1\star f_2 = \sum_{f_1=f_1' f_2'' \atop f_2=f_1' f_2''} 
(f_1'\star f_2') M(f_1'';f_2'')
$$
contains exactly all possible shuffles mixed with all possible
graftings. This is guaranteed by the inductive hypothesis on the
factors $f_1'\star f_2'$ and $M(f_1'';f_2'')$, and the fact that the
sum runs over all decompositions of $f_1$ and $f_2$ into two factors, 
including the trivial ones. 
\end{proof}
  
Here are some examples of graftings and shuffles: 
\begin{align*}
M\left(\pstree{\tbc{b}{x}}{\tbc{a}{y}} ; \pstree{\tbc{b}{u}}{} \hspace{.3cm}
  \pstree{\tbc{b}{v}}{} \right)
&= \pstree{\tbc{b}{x}}{\tbc{a}{u}\tbc{a}{v}\tbc{a}{y}}
+\pstree{\tbc{b}{x}}{\tbc{a}{u}\pstree{\tbc{r}{y}}{\tbc{a}{v}}}
+\pstree{\tbc{b}{x}}{\tbc{a}{u}\tbc{a}{y}\tbc{a}{v}} +
\pstree{\tbc{b}{x}}{\pstree{\tbc{r}{y}}{\tbc{a}{u}\tbc{a}{v}}} +
\pstree{\tbc{b}{x}}{\pstree{\tbc{l}{y}}{\tbc{a}{u}}\tbc{a}{v}}+
\pstree{\tbc{b}{x}}{\tbc{a}{y}\tbc{a}{u}\tbc{a}{v}} 
\\ 
\pstree{\tbc{b}{x}}{} * \left ( \pstree{\tbc{b}{u}}{} \hspace{.3cm}
  \pstree{\tbc{b}{v}}{\tbc{a}{w}} \right ) 
&= \pstree{\tbc{b}{x}}{} \hspace{.3cm}\pstree{\tbc{b}{u}}{} \hspace{.3cm}
  \pstree{\tbc{b}{v}}{\tbc{a}{w}}+
  \pstree{\tbc{b}{u}}{}\hspace{.3cm} \pstree{\tbc{b}{x}}{}
  \hspace{.3cm} \pstree{\tbc{b}{v}}{\tbc{a}{w}} +\pstree{\tbc{b}{u}}{}
  \hspace{.3cm} \pstree{\tbc{b}{v}}{\tbc{a}{w}}
  \hspace{.3cm}\pstree{\tbc{b}{x}}{}
  +\pstree{\tbc{b}{x}}{\tbc{a}{u}}\hspace{.3cm}
  \pstree{\tbc{b}{v}}{\tbc{a}{w}} + \pstree{\tbc{b}{x}}{\tbc{a}{u}
    \pstree{\tbc{r}{v}}{\tbc{a}{w}}} + \pstree{\tbc{b}{u}}{}
  \hspace{.3cm}\pstree{\tbc{b}{x}}{\pstree{\tbc{r}{v}}{\tbc{a}{w}}}. 
\end{align*}

According to the results of Loday and Ronco, in particular paragraph
3.14 in \cite{LodayRonco2008}, we can conclude that 

\begin{corollary} 
The brace algebra $R^{\calT}$, endowed with $M$, is the free brace 
algebra over $R$. 
\end{corollary}

\begin{remark}
The algebras $R^{\calT}$ and $T(R^{\calT})$, or their duals, or their
non-planar and commutative versions, appeared in really many recent
works. We quote, but the list is surely not exhaustive, the works 
by R.~Grossman and R.~G.~Larson \cite{GrossmanLarson}, 
F.~Chapoton \cite{Chapoton}, L.~Foissy \cite{Foissy}, D.~Guin et
J.-M.~Oudom \cite{GuinOudom}, M.~Ronco \cite{Ronco} and 
J.-L.~Loday and M.~Ronco \cite{LodayRonco2008}. 
\end{remark}
\bigskip 

Now suppose that $R$ itself is a brace algebra, and therefore $T(R)$
is a right-sided CHA. We are going to show that the $\star$ product
defined in $T(R)$ can be found by computing its lifting to
$T(R^{\calT})$. 

Let $\iota$ be the linear morphism from $T(R)$ to $T(R^{\calT})$ defined
by $\iota(1)=1$ and 
$$
\iota(x_1\cdots x_n)=B_+^{x_1}(1) \cdots B_+^{x_n}(1)=
\pstree{\tbc{b}{x_1}}{}\cdots \pstree{\tbc{b}{x_n}}{} . 
$$
It is clear that $\iota$ preserves the concatenation in $T(R)$ and
$T(R^{\calT})$, and that it identifies $R$ with $R^{\bullet}\subset
R^{\calT}$. 

Note that the map $\iota$ does not preserve the brace product, in fact 
$$
M\left ( \pstree{\tbc{b}{x}}{} ; \pstree{\tbc{b}{y_1}}{} \dots 
\pstree{\tbc{b}{y_q}}{} \right ) 
= \pstree{\tbc{b}{x}}{\tbc{l}{y_1} \Tfan[linestyle=dotted]\tbc{r}{y_q}} 
\neq \iota(M(x;y_1\cdots y_q))=\pstree{\tbc{b}{M(x;y_1\cdots y_q)}}{}, 
$$
and therefore it does not preserve the product $\star$ in $T(R)$ and
$T(R^{\calT})$. 
We show that $\iota$ has an inverse map which allows to compute the
product $\star$ in $T(R)$ starting from that in $T(R^{\calT})$. 

\begin{theorem}
\label{thm-ouf}
Let $\mu$ the linear morphism from $T(R^{\calT})$ to $T(R)$ defined recursively by
\begin{enumerate}
\item 
$\mu(1)=1$,
\item 
$\mu(t_1 \ldots t_n)=\mu(t_1)\ldots \mu(t_n)$,
\item 
$\mu\big(B_+^x(s_1\ldots s_q)\big)=M_{1q}\big(x;\mu(s_1\ldots s_q)\big)$.
\end{enumerate}
Then, for any $X,Y \in T(R)$, we have 
\begin{align}
\label{formula-ouf}
X\star Y = \mu\big(\iota(X)\star\iota(Y)\big). 
\end{align}
\end{theorem}

\begin{proof}
The proof is straightforward and relies on the recursive formulas 
in brace algebras and on the fact that $\mu \circ \iota=Id_{T(R)}$. 
If $X=x_1 \cdots x_p$ and $Y=y_1 \cdots y_q$ then
\begin{align*}
X \star Y & = \sum_{X=X^1 X^2 \atop Y=Y^1 Y^2} (X^1\star Y^1) M(X^2;Y^2)
\end{align*}
Since $R$ is a brace algebra, either $X^2=1$ or $X^2=x_p$. 
\begin{itemize}
\item 
If $X^2=1$, then
$$
M(1,Y^2)=Y^2=\mu(\iota(Y^2))=\mu\big(M(1;\iota(Y^2)\big). 
$$
\item 
If $X^2=x_p$, then
$$
M(x_p,Y^2)=\mu(B_+^{x_p}(\iota(Y^2)))=\mu\big(M(\iota(x_p);\iota(Y^2))\big). 
$$
\end{itemize}
Recursively, we get
\begin{align*}
X \star Y &=\sum_{X=X^1 X^2 \atop Y = Y^1 Y^2} (X^1 \star Y^1) M(X^2;Y^2)\\
&=\sum_{X=X^1 X^2 \atop Y=Y^1 Y^2} \mu\big(\iota(X^1)\star\iota(Y^1)\big) 
\mu\big(M(\iota(X^2);\iota(Y^2))\big) \\
&= \mu\left(\sum_{\iota(X)=\iota(X^1)\iota(X^2) \atop 
\iota(Y)=\iota(Y^1)\iota(Y^2)} 
\big(\iota(X^1)\star\iota(Y^1)\big) M\big(\iota(X^2);\iota(Y^2)\big)\right)\\
&= \mu\big(\iota(X)\star\iota(Y)\big).  
\end{align*}
\end{proof}

This result means that one can compute the $\star$ product in $T(R)$
with the help of the $\star$ product in $T(R^{\calT})$. This last product,
given by shuffle and grafting of trees, is quite easy to handle. 
Let us illustrate the application with an example:
\begin{align*}
x \star (yz) &= \mu \left ( \pstree{\tbc{b}{x}}{}\star 
\left(\pstree{\tbc{b}{y}}{} \hspace{.3cm} 
\pstree{\tbc{b}{z}}{}\right)\right)\\
&= \mu \left (
  \pstree{\tbc{b}{x}}{} \hspace{.3cm} \pstree{\tbc{b}{y}}{}
  \hspace{.3cm}  \pstree{\tbc{b}{z}}{}
  +\pstree{\tbc{b}{y}}{} \hspace{.3cm} \pstree{\tbc{b}{x}}{}
  \hspace{.3cm}  \pstree{\tbc{b}{z}}{}
  +\pstree{\tbc{b}{y}}{} \hspace{.3cm} \pstree{\tbc{b}{z}}{}
  \hspace{.3cm}  \pstree{\tbc{b}{x}}{} 
  +\pstree{\tbc{b}{x}}{\tbc{a}{y}} \hspace{.3cm} \pstree{\tbc{b}{z}}{}+
  \pstree{\tbc{b}{y}}{} \hspace{.3cm} \pstree{\tbc{b}{x}}{\tbc{a}{z}} +
  \pstree{\tbc{b}{x}}{\tbc{l}{y} \tbc{r}{z}}  \right ) \\
&= xyz+ yxz+yzx+M_{11}(x;y)z+ yM_{11}(x;z)+ M_{12}(x;yz). 
\end{align*}
This computation can be compared with that made in section~\ref{rs-CHA}.

\begin{remark}
The result of theorem (\ref{thm-ouf}) can be obtained also as a consequence 
of the equivalence between the categories of brace and dendriform algebras 
proved by Chapoton in \cite{Chapoton}. Since $\mu$ is by assumption a brace 
homomorphism between $R^{\calT}$ and $R$, this equivalence allows to state that 
for any $f_1$ and $f_2$ in $T(R^{\calT})$ we have 
$\mu(f_1)\star\mu(f_2)=\mu(f_1\star f_2)$. If we then choose $f_1=\iota(X)$ 
and $f_2=\iota(Y)$, we obtail formula (\ref{formula-ouf}). 
\end{remark}

\begin{remark}
The algebra $T(R^{\calT})$ is the dual of Foissy's Hopf algebra of decorated
plane trees, cf.~\cite{Foissy}. Foissy's Hopf algebra on rooted trees has a
``universal property'', in the sense that it is endowed with the
Hochschild cocycle operator $B_+$. Note that the property that we show
for $T(R^{\calT})$ is not the dual property of that one.  
\end{remark}

%%%%%%%%%%%%%%%%%%%%%%%%%%

\subsection{Right-sided CHAs from dual Hopf algebras}

In the next sections we determine some right-sided CHAs $T(R)$ from the 
dualization of given Hopf algebras $T(V)$ where the vector space $V$
is itself graded. 
In this case, the vector space $R$ coincides with the dual space
$V^*$, and we give a general tool to compute the brace structure on $R$ 
starting from the coproduct on $T(V)$. This tool is then used in the examples 
presented in the rest of the paper.  

Let $\H$ be a graded Hopf algebra that is free as an algebra, 
finite dimensional in all degrees, and with generators 
which are themselves graded. 
If we denote by $V$ be the vector space spanned by all the
generators of $\H$ (infinite but countable many), the Hopf algebra 
can be given as $\H=T(V)=\bigoplus_{n=0}^\infty V^{\otimes n}$.

The dual Hopf algebra $\H^*$ is defined as follows.
As a vector space, it is the graded linear dual vector space of
$\H$, that is $\H^*=T(V^*)$, where $V^*$ is the graded linear dual 
vector space of $V$. 
If we denote by $\{v_n, n\geq 1\}$ a generic basis of $V$, that is a
set of generators of $\H$, the element $v_n^*$ dual to each $v_n$ 
is defined by the pairing $\langle v_n^*| v_m\rangle=\delta_{n,m}$.
Therefore, the dual of a generic basis element $v_{n_1}\cdots v_{n_i}$ in
$\H$ is the element $v^*_{n_1}\cdots v^*_{n_i}$ in $\H^*$ 
defined by the pairing
\begin{eqnarray*}
\langle v^*_{n_1}\cdots v^*_{n_i}|v_{m_1}\cdots v_{m_j}\rangle 
&=& \delta_{i,j} 
\langle v^*_{n_1}|v_{m_1}\rangle \dots
\langle v^*_{n_i}|v_{m_i}\rangle.
\end{eqnarray*}

As a coalgebra, $\H^*$ inherits the coproduct $\D^*$ dual to the product
of $\H$, that is, such that 
\begin{align*} 
\big\langle \D^*(v^*_{n_1}\cdots v^*_{n_i})| 
(v_{m_1}\cdots v_{m_j}) \otimes (v_{p_1}\cdots v_{p_k}) \big\rangle 
&= 
\big\langle v^*_{n_1}\cdots v^*_{n_i}| 
v_{m_1}\cdots v_{m_j} v_{p_1}\cdots v_{p_k} \big\rangle. 
\end{align*}
Since $\H$ is free and unital as an algebra, $\H^*$ is cofree and
counital as a coalgebra, that is, $\D^*$ is the deconcatenation 
coproduct 
\begin{align*} 
\D^*(v^*_n) &= v^*_n \otimes 1 + 1 \otimes v^*_n, \\ 
\D^*(v^*_{n_1}\cdots v^*_{n_k}) &= \sum_{i=0}^k 
v^*_{n_1}\cdots v^*_{n_i}\otimes v^*_{n_{i+1}}\cdots v^*_{n_k}, 
\end{align*}
and with the counit $\varepsilon = 1^*$. 
In this expression, $1$ is not (yet) a unit, it denotes the generator of
$(V^*)^{\otimes 0}$. 
The set of primitive elements of $\H^*$ is then obviously the dual
$V^*$ of the set of generators of $\H$.

Finally, as an algebra, $\H^*$ inherits the product $\star$ dual to the
coproduct $\D$ of $\H$, that is, the product such that 
\begin{align*} 
\big\langle (v^*_{n_1}\cdots v^*_{n_i})\star (v^*_{m_1}\cdots v^*_{m_j}), 
v_{p_1}\cdots v_{p_k} \big\rangle 
&= \big\langle v^*_{n_1}\cdots v^*_{n_i}\otimes\ v^*_{m_1}\cdots v^*_{m_j}, 
\D(v_{p_1}\cdots v_{p_k}) \big\rangle. 
\end{align*}

The right-sided condition (\ref{right-sided-cofree}) for the 
Hopf algebra $T(V^*)$ is equivalent to require that the coproduct 
$\D$ on $T(V)$ induces a right coaction of the Hopf algebra $T(V)$ 
on the subspace $T^{\leq p}(V)=\oplus_{n\leq p} V^{\otimes n}$, that is 
\begin{align}
\label{right-sided-free}
\D\big(T^{\leq p}(V)\big) \quad\subset\quad T^{\leq p}(V) \otimes T(V). 
\end{align}
Explicitly, this means that for any element 
$a=v_{n_1}\cdots v_{n_p}\in V^{\otimes p}$, the coproduct 
$\D(a)= \sum a_{(1)} \otimes a_{(2)}$ (in Sweedler's notation) 
produces terms with $a_{(1)} \in T^{\leq p}(V)$. 
Since the coproduct is an algebra morphism, and the product is free, 
it is sufficient that this property holds for $p=1$: for any generator
$v\in V$, the coproduct $\D(v)= \sum v_{(1)} \otimes v_{(2)}$ 
produces terms with $v_{(1)} \in T^{\leq 1}(V)= K \oplus V$ 
(cf. \cite{LodayRonco2008}, Section 3.6). 

Finally, the brace product on $V^*$ is the projection of $\star$ onto $V^*$. 
The projection $\pi:T(V^*)\longrightarrow V^*$, in this case, is the map 
\begin{align*} 
\pi(a^*) &= \sum_m \langle a^* | v_m \rangle\ v_m^*, 
\end{align*}
for any $a^*\in T(V^*)$. 
Therefore, the maps $M_{1q}:V^*\otimes (V^*)^{\otimes q}\longrightarrow V^*$ 
are simply given by 
\begin{align}
M_{1q}(v^*_n;a^*) &=
\sum_m \langle v^*_n \star a^* | v_m \rangle\ v^*_m =
\sum_m \langle v^*_n \otimes a^* | \D v_m \rangle\ v^*_m.
\label{M1q}
\end{align}
This sum if finite. In fact, let us denote by $|a|$ [resp. $|a^*|$] the degree 
of the elements in the graded Hopf algebra $T(V)$ [resp. $T(V^*)$]. 
We recall that if $|v_n|$ denotes the degree of a generic basis element of
$V$, then the degree of an element $v_{n_1}\cdots v_{n_p}\in T(V)$
is given by $|v_{n_1}\cdots v_{n_p}|=|v_{n_1}|+\cdots +|v_{n_p}|$. 
Then the sum over $m$ in Formula~(\ref{M1q}) is limited by the fact that 
$|v^*_m|=|v_m|=|v_n^*|+|a^*|$.

%%%%%%%%%%%%%%%%%%%%%%%%%%%%%%%%%%%%%%%%%%%%%%%%%%%%%%%%%%%%%%%%%%%%

\section{The right-sided combinatorial structure 
of the Hopf algebra $(\Hdnc)^*$}

Let $\Hdnc$ be the non-commutative Hopf algebra of formal 
diffeomorphisms, as defined in \cite{VanDerLaan} or \cite{BFK} 
(where it was denoted by $\Hd$). 
It is the graded and connected Hopf algebra 
$\Q\langle v_1,v_2,...\rangle=T(V)$ on the generators 
$V=\Span\{v_1,v_2,...\}$ graded by $|v_n|=n$, 
considered with the free (tensor) product, 
and endowed with the coproduct 
\begin{align}
\label{Hdif-coproduct}
\Dd(v_n) &= \sum_{m=0}^n v_m \otimes 
\underset{k_0+k_1+\dots+k_m=n-m}{\sum_{k_0,k_1,\dots,k_m\geq 0}} 
v_{k_0} v_{k_1}\dots v_{k_m},  
\end{align}
where we set $v_0=1$, and with the counit $\varepsilon(v_n)=\delta_{n,0}$. 
This Hopf algebra is a free-associative CHA, and the formula 
(\ref{Hdif-coproduct}) clearly says that the coproduct $\Dd$ 
satisfies the right-sided condition (\ref{right-sided-free}). 

Let $(\Hdnc)^*=T(V^*)$ be its dual Hopf algebra, as described in section 
\ref{section1}, with primitive elements in $V^*=\Span\{v^*_1,v^*_2,...\}$. 
Then $(\Hdnc)^*=T(V^*)$ is a right-sided cofree-coassociative CHA, 
and moreover a dendriform algebra, and $V^*$ is a brace algebra. 

\begin{proposition}
The brace product on $V^*$ which induces the CHA structure on
$(\Hdnc)^*$ is given by 
\begin{align*}
\{v_{n}^*;v_{m_1}^*\cdots v_{m_q}^*\} 
= M_{1q}(v_{n}^*;v_{m_1}^*\cdots v_{m_q}^*) 
&= \binom{n+1}{q}\ v_{n+m_1+\cdots +m_q}^*. 
\end{align*}
\end{proposition}

\begin{proof}
According to eq.~(\ref{M1q}), for any
$v_n\in V$ and any $a\in \Hdnc$, we have
\begin{align*}
M_{1q}(v_n^*;a^*) &= 
\sum_{m=1}^\infty \langle v_n^* \otimes  a^*\ |\ \Dd v_m\rangle\ v_m^*.
\end{align*}
If we put $v_0=1$, the coproduct of $v_m$ becomes~\cite{BFK}
\begin{align*}
\Dd v_m &= \sum_{k=0}^m v_k \otimes Q^{(k)}_{m-k}(v),
\end{align*}
with
\begin{align*}
Q^{(k)}_{m-k}(v) &= \sum_{l=1}^{k+1} \binom{k+1}{l}
\underset{j_1+\dots+j_l=m-k}{\sum_{j_1,\dots,j_l > 0}} 
v_{j_1}\dots v_{j_l}.
\end{align*}
Therefore, 
\begin{align*}
M_{1q}(v_n^*;a^*) &= 
\sum_{m=1}^\infty \langle a^*\ |\ Q^{(n)}_{m-n}(v) \rangle\ v_m^*.
\end{align*}
We take $a=v_{m_1}\dots v_{m_q}$ and we evaluate
$\langle v_{m_1}^*\dots v_{m_q}^*\ |\ Q^{(n)}_{m-n}(v) \rangle$.
From the definition of $Q^{(n)}_{m-n}(v)$, we must
have $l=q$, $j_1=m_1,\dots,j_q=m_q$.
Therefore, $m-n=j_1+\dots+j_q=m_1+\dots+m_q$ and 
\begin{align*}
M_{1q}(v_n^*;v_{m_1}^*\dots v_{m_q}^*) &= 
\binom{n+1}{q}\ v_{n+m_1+\dots+m_q}^*.
\end{align*}
\end{proof}

%%%%%%%%%%%%%%%%%%%%%%%%%%%%%%%%%%%%%%%%%%%%%%%%%%%%%%%%%%%%%%%%%%%%

\section{The right-sided combinatorial structure 
of the Hopf algebra $(\Hanc)^*$}

Let $Y$ be the set of planar binary rooted trees, that is, planar graphs 
without loops and a preferred external edge called the root. For instance, 
\begin{align*}
\treeO, \quad \treeA, \quad \treeAB, \quad \treeBA, \quad 
\treeABC, \quad \treeBAC, \quad \treeACA, \quad \treeCAB, \quad \treeCBA. 
\end{align*}
For a tree $t$, we denote by $|t|$ the number of its internal vertices. 
We denote by $Y_n$ the set of planar binary trees with $n$ internal 
vertices, so that $Y=\bigcup_{n=0}^\infty Y_n$. 

Let $\Hanc$ be the non-commutative lift of the Hopf algebra $\Ha$ introduced 
in \cite{BFqedren,BFqedtree} to describe the renormalization of the 
electric charge in the perturbative expansion of quantum electrodynamics 
based on planar binary trees. 
It is the graded and connected Hopf algebra $\Q Y$ spanned by all planar 
binary trees, endowed with the product {\em over}, 
cf.~\cite{Loday2002}: 
given two planar binary trees $s,t\neq\treeO$, the tree $s$ over $t$, 
denoted by $s\over t$, is the tree obtained by grafting the root of $s$ 
over the left-most leaf of $t$, that is, 
\begin{align*}
s \over t &= \lgraft{t}{s}. 
\end{align*}
This product is not commutative, and the root tree $\treeO$ is the unit. 
Moreover, any tree can be decomposed as the over product of subtrees which 
have nothing branched on their left-most leaf. 
If we set $v(t)=\rgraft{\treeA}{t}$ to denote these trees, the algebra 
$\Hanc$ is in fact isomorphic to the free algebra $T(V)$, where 
$V=\Span\{v(t),t\in Y\}$. 

The coproduct $\Da:\Hanc \longrightarrow \Hanc \otimes\Hanc$ can be
described in an elegant way in the form proposed by P.~Palacios in her
Master Thesis \cite{Palacios}. 
For this, we need some notation. Let us denote by 
$\vee :Y_n \times Y_m \longrightarrow Y_{n+m+1}$ the operation 
which grafts two trees on a new root, that is
\begin{align*}
s \vee t &= \lrgraft{s}{t}, 
\end{align*}
and call dressed comb the tree $\gamma(t_1,\dots,t_k)$ recursively
given by 
\begin{align*}
\gamma(t) &= t\vee \treeO = \lgraft{\treeA}{t},\\
\gamma(t_1,t_2,\dots,t_k) &= t_1\vee \gamma(t_2,\dots,t_k) 
= \combRgraft{t_1}{t_2}{t_k}.
\end{align*} 
Then, any tree $t\neq\treeO$ can be written as 
$t=\gamma(t_1,t_2,\dots,t_k)$ for some suitable trees $t_i$. 

The coproduct $\Da$ is the unital algebra homomorphism defined on 
the generators as 
\begin{align}
\nonumber 
\Da \treeA &= \treeO \otimes \treeA + \treeA \otimes \treeO, \\ 
\label{Da-def}
\Da v(t) &= \treeO \otimes  v(t)
+ \sum\ v(\gamma(t_1\i1,t_2\i1,\dots,t_k\i1)) 
\otimes t_1\i2\over t_2\i2\over \cdots\over t_k\i2, 
\end{align}
where $t=\gamma(t_1,t_2,\dots,t_k)$ and 
$\sum t_i\i1 \otimes t_i\i2 = \Da(t_i)$ is the Sweedler notation 
for the coproduct applied to each subtree $t_i$. 

Finally, the counit $\varepsilon:\Hanc\longrightarrow \Q$ is the unital 
algebra morphism with value $\varepsilon(v(t))=0$ on the generators. 

Let $(\Hanc)^*=T(V^*)$ be its dual Hopf algebra, with primitive elements 
given by $V^*=\Span\{v(t^*),t\in Y\}$, where $t^*$ is the dual form of the 
tree $t$. Then $(\Hanc)^*$ is a right-sided cofree-coassociatibe CHA, and 
$V^*$ is a brace algebra. We describe here its brace structure. 

\begin{theorem}
The brace product $M_{1q}$ on $V^*$ which induces the CHA structure on
$(\Hdnc)^*$ is given for any $q>0$ by
\begin{align*}
M_{1q}(\treeA^*;v(s_1)^*\cdots v(s_q)^*) &= 0 
\end{align*}
and for any $t\neq \treeO$ by 
\begin{align*}
M_{1q}(v(t)^*;v(s_1)^*\cdots v(s_q)^*) &= 
\sum\ v\big(\gamma(p^*\i1,\dots,p^*\i{k})\big)^*, 
\end{align*}
where for $t=\gamma(t_1,\dots,t_k)$ we set 
$p^*=(t_1^*\cdots t_l^*)\star (v(s_1)^*\cdots v(s_q)^*)$, 
and where we use the Sweedler notation 
\begin{align*}
\sum\ p^*\i1\otimes p^*\i2\otimes\cdots\otimes p^*\i{k}
&=(\D^*)^{k-1}\ p^*
\end{align*}
for the deconcatenation coproduct $\D^*$ applied $k-1$ times to $p^*$. 
\end{theorem}

\begin{proof}
According to Eq.~(\ref{M1q}), if we denote $a^*=v(s_1)^*\cdots v(s_q)^*$, 
the brace product on $V^*$ can be computed as 
\begin{align}
\nonumber
M_{1q}\big(v(t)^*;a^*\big) &=
\sum_{u\in Y} \big\langle v(t)^*\star a^* | v(u)\big\rangle\ v(u)^* 
= \sum_{u\in Y} \big\langle v(t)^*\otimes a^* | \Da v(u)\big\rangle\
v(u)^* \\ 
&= \big\langle v(t)^*\otimes a^* | \Da v(\treeO)\big\rangle\ v(\treeO)^* 
+ \sum_{u\neq\treeO} 
\big\langle v(t)^*\otimes a^* | \Da v(u)\big\rangle\ v(u)^* . 
\label{term2}
\end{align}
The first term gives 
\begin{align*}
\big\langle v(t)^*\otimes a^* | \treeA\otimes\treeO
+\treeO\otimes\treeA \big\rangle\ \treeA^* &= 0, 
\end{align*}
because $\langle a^*\ |\ \treeO\rangle =0$ for $q>0$, 
and $\langle v(t)^* |\ \treeO \rangle = 0$ for any $t$. 

Then we suppose that $u\neq\treeO$ and we evaluate the second term 
of (\ref{term2}) using the expression (\ref{Da-def}) for $\Da v(u)$. 
The first part gives  
\begin{align*}
\big\langle v(t)^*\otimes a^* |\ \treeO \otimes v(u) \big\rangle &=0 
\end{align*} 
for any choice of $t$. 

The second part gives
\begin{align} 
\nonumber
&\big\langle v(t)^*\otimes a^* | v(\gamma(u_1\i1,u_2\i1,\dots,u_l\i1)) 
\otimes u_1\i2\over u_2\i2 \over\cdots\over u_l\i2 \big\rangle \\ 
&\hspace{2cm} = 
\big\langle v(t)^* | v(\gamma(u_1\i1,u_2\i1,\dots,u_l\i1))\big\rangle
\big\langle a^*| u_1\i2\over u_2\i2 \over\cdots\over u_l\i2 \big\rangle ,  
\label{M1q-2}
\end{align}
where we write $u$ as $\gamma(u_1,...,u_l)$. 
The only primitive element $v(t)$ which can not be of the form 
$v\big(\gamma(u_1\i1,u_2\i1,\dots,u_k\i1)\big)$ is $v(\treeO)=\treeA$. 
In this case, the brace product is obviously 
\begin{align*}
M_{1q}\big(\treeA^*;a^*\big) &= 0 
\end{align*}
for any $q>0$. 
Finally let us consider the case $t\neq\treeO$. 
Then $t = \gamma(t_1,\dots,t_k)$,
and the first pairing of (\ref{M1q-2}) is non-zero only if 
\begin{align*}
\gamma(u_1\i1,u_2\i1,\dots,u_l\i1) &= t = \gamma(t_1,\dots,t_k), 
\end{align*}
that is, if and only if $l=k$ and 
\begin{align*}
u_1\i1=t_1\ , \quad u_2\i1=t_2\ ,\quad\dots\ ,\quad u_k\i1=t_k\ .  
\end{align*}
Therefore we have
\begin{align*}
M_{1q}(v(t)^*;a^*) &= \sum_{u_1,\dots,u_k}
\big\langle v(t)^*\ |\ v\big(\gamma(u_1\i1,\dots,u_k\i1)\big) \big\rangle
\big\langle a^*\ |\ u_1\i2\over\cdots\over u_l\i2\big\rangle 
v\big(\gamma(u_1,\dots,u_k)\big)^* \\ 
&= \sum_{u_1,\dots,u_k}
\big\langle t_1^*\cdots t_k^*\ |\ u_1\i1\cdots u_k\i1 \big\rangle 
\big\langle a^*\ |\ u_1\i2\over\cdots\over u_l\i2\big\rangle 
v\big(\gamma(u_1,\dots,u_k)\big)^* \\ 
&= \sum_{u_1,\dots,u_k}
\big\langle t_1^*\cdots t_k^* \otimes a^*\  |\  
\Da(u_1\over\cdots\over u_k) \big\rangle 
v\big(\gamma(u_1,\dots,u_k)\big)^* \\ 
&= \sum_{u_1,\dots,u_k}
\big\langle (t_1^*\cdots t_k^*) \star a^*\  |\ 
u_1\over\cdots\over u_k \big\rangle 
v\big(\gamma(u_1,\dots,u_k)\big)^* .  
\end{align*}
Now, note that the deconcatenation coproduct $\D^*$ in $(\Hanc)^*$ 
is dual to the product $\over$ in $\Hanc$, that is 
\begin{align*}
\big\langle a^*\ | b\over c \big\rangle &= 
\big\langle \D^*(a^*)\ | b \otimes c \big\rangle. 
\end{align*}
Let us call $p^*= (t_1^*\cdots t_k^*) \star a^*$. 
If we apply $k-1$ times the deconcatenation coproduct $\D^*$ to $p^*$,
and we use the Sweedler notation 
$\sum p^*\i1\otimes p^*\i2\otimes\cdots\otimes p^*\i{l}=(\D^*)^{k-1}\ p^*$, 
we obtain the final result 
\begin{align*}
M_{1q}(v(t)^*;a^*) 
&= \sum_{u_1,\dots,u_k} \sum\ 
\big\langle p^*\i1\otimes p^*\i2\otimes\cdots\otimes p^*\i{k}\ |\ 
u_1\otimes\cdots\otimes u_k \big\rangle 
v\big(\gamma(u_1,\dots,u_k)\big)^* \\ 
&= \sum\ v\big(\gamma(p^*\i1,\dots,p^*\i{k})\big)^* .  
\end{align*}
\end{proof}

\begin{examples}
The simplest example is given by $M_{11}({\treeBA}^*;{\treeA}^*)$.
We have $\treeBA=v(\treeA)=v(\gamma(\treeO))$, so that
$k=1$ and $t_1=\treeO$. 
Therefore, $p^*=\treeO^*\star \treeA^*=\treeA^*$ and 
\begin{align*}
M_{11}(\treeBA^*;\treeA^*) &= v(\gamma(\treeA^*))=v(\treeAB^*)=\treeCAB^*.
\end{align*}

Similarly
\begin{align*}
M_{11}(\treeBA^*;\treeBA^*) &= v(\gamma(\treeBA^*))=v(\treeBAC^*)=\treeDBAC^*.
\end{align*}

Another simple example is given by $M_{11}(\treeCAB^*;\treeA^*)$.
We have $\treeCAB= v(\gamma(\treeA))$, so that $k=1$ and $t_1=\treeA$. 
Since for any $t$ formula (\ref{M1q}) gives 
\begin{align*}
\treeA^* \star v(t)^* &= 
v(t)^*\treeA^* + \treeA^* v(t)^* + M_{1q}(\treeA^*;v(t)^*) 
= v(t)^*\treeA^* + \treeA^* v(t)^*, 
\end{align*}
because $M_{1q}(\treeA^*;v(t)^*)=0$, we get  
$p^*=\treeA^*\star\treeA^*= 2 \treeAB^*$, and finally 
\begin{align*}
M_{11}(\treeCAB^*;\treeA^*) &= 2 v(\gamma(\treeAB))^*
=2v(\treeABC)^*=2\treeDABC^*.
\end{align*}

Now let us compute $M_{11}(\treeCBA^*;\treeA^*)$.
We have $\treeCBA=v(\gamma(\treeO,\treeO))$ so that $k=2$ and 
$t_1=t_2=\treeO$. 
Since $\treeO^*\treeO^*=\treeO^*$, we have 
$p^*=\treeO^*\star\treeA^*=\treeA^*$ and therefore 
\begin{align*}
M_{11}(\treeCBA^*;\treeA^*) &=
v(\gamma(\treeA,\treeO))^*+v(\gamma(\treeO,\treeA))^*  
= v(\treeACA)^*+v(\treeCAB)^* = \treeDACA^* + \treeDCAB^*.
\end{align*}

Finally, let us compute an example of $M_{1q}$ with $q>1$. 
We consider $M_{12}(\treeBA^*;\treeA^*\treeA^*)$. 
We have $\treeBA=v(\gamma(\treeO))$, so that $k=1$ and $t_1=\treeO$. 
Therefore, $p^*=\treeO^*\star (\treeA^* \treeA^*)=\treeAB^*$ and 
\begin{align*}
M_{12}(\treeBA^*;\treeA^*\treeA^*) &=
v(\gamma(\treeAB))^*=v(\treeCAB)^* = \treeDCAB^*.
\end{align*}
Similarly
\begin{align*}
M_{12}(\treeCBA^*;\treeA^*\treeA^*) &=
v(\gamma(\treeAB,\treeO))^*+v(\gamma(\treeA,\treeA))^*
+v(\gamma(\treeO,\treeAB))^* \\ 
& = v(\treeABDA)^* + v(\treeADAB)^* + v(\treeDABC)^* 
= \treeEABDA^* + \treeEADAB^* + \treeEDABC^*.
\end{align*}
\end{examples}

\begin{remark}
The brace structure can be described as follows: 
\begin{align*}
M_{1q}(v(t)^*;v(s_1)^*\cdots v(s_q)^*) &= 
\sum_u\ v(u)^*, 
\end{align*}
where the sum runs over any tree $u$ obtained by branching the trees
$v(s_i)$ or any block $v(s_i)\over\cdots\over v(s_j)$ on the 
$\backslash$-leaves of $t$ or inside the $\backslash$-branches of
$t$, by preserving the order from left to right of the trees $v(s_i)$
and of the $\backslash$-branches of $t$. 
\end{remark}

%%%%%%%%%%%%%%%%%%%%%%%%%%%%%%%%%%%%%%%%%%%%%%%%%%%%%%%%%%%%%%%%%%%%

\section{The right-sided combinatorial structure of the Hopf algebra 
$T(\T(B))^*$}

In the paper \cite{BrouderSchmitt}, C.~Brouder and W.~Schmitt
introduced a non-commutative version $\H$ of a renormalization Hopf algebra
inspired by Pinter's work on the Epstein-Glaser renormalization of
quantum fields on the configuration space. 
In this section we describe the brace structure on the set of
primitive elements of its dual Hopf algebra $\H^*$. 

As an algebra, $\H$ is the tensor algebra $T(V)$ on a set of generators 
$V=\T(B) =\oplus_{n\geq 1}\ B^{\otimes n}$ which is the augmentation
ideal of the tensor algebra over a given bialgebra $B$. 

In order to apply the dualization procedure described in 
Section~\ref{section1}, we have to assume that the bialgebra $B$
admits at most a countable basis, and that this has been fixed. 
Let us denote by $\{x_n\}$ the chosen basis of $B$. 
We also denote the product (in $B$) by a dot, namely $x\cdot y$ for 
any $x,y\in B$, and the coproduct (in $B$) through the Sweedler 
notation $\d(x)=\sum\ x\i1 \otimes x\i2$. 

In $V=\T(B)$, we denote the tensor products by tuples. Therefore the 
elements $v_{n_1,...,n_k}=(x_{n_1},\dots,x_{n_k})$ form a basis of $V$. 
On $V$ we consider two coproducts. 
The first one is the extention of the coproduct on $B$ as an algebra 
morphism, namely 
\begin{align*}
\d(x,y,\dots,z) &= \sum\ (x\i1,y\i1,\dots,z\i1) 
\otimes (x\i2,y\i2,\dots,z\i2). 
\end{align*}
For a generic element $v$ of $V$, we use the same Sweedler notation 
$\d(v)=\sum\ v\i1\otimes v\i2$ as on $B$. 
The second coproduct is the (reduced) deconcatenation coproduct 
\begin{align*}
\D(x,y,\dots,z) &= (x) \otimes (y,\dots,z) + (x,y) \otimes (\dots,z)
+ \cdots + (x,y,\dots) \otimes (z). 
\end{align*}
For this coproduct we use a modified Sweedler notation 
$\D(v)=\sum\ v\j1\otimes v\j2$. 

Finally, in $\H=T(V)$, we denote the tensor product by the
concatenation. A generic element in $\H$ is then of the form 

$a=uv\cdots w$, where $u,v,...,w\in V$, and for any elements 
$a,b$ we denote their product (in $T(V)$) by $ab$. 
Both coproducts $\d$ and $\D$ can be extended to $\H$. 
The renormalization coproduct is then the algebra homomorphism defined
on a generator $v\in V$ as 
\begin{align}
\label{DH-def}
\DH(v) &= \sum_{n=1}^\infty \sum\ 
\big(\mu(v\j1\i1),\dots,\mu(v\j{n}\i1)\big) 
\otimes v\j1\i2 \cdots v\j{n}\i2, 
\end{align}
where $\sum\ v\j1\otimes\cdots\otimes v\j{n}=\D^{n-1}(v)$ is the 
result of the deconcatenation applied $n-1$ times to $v$, 
and $\mu:\T(B)\longrightarrow B$, 
$\mu(x,y,...,z)=x\cdot y\cdot\ldots\cdot z$ is the product of $B$
extended to many factors. 

Note that for any $v\in V$, the left tensor factor 
$\big(\mu(v\j1\i2),\dots,\mu(v\j{n}\i2)\big)$ of (\ref{DH-def}) 
belongs to $V$ (in fact it belongs to $B^{\otimes n}$), 
while the right factor $v\j1\i1 \cdots v\j{n}\i1$ belongs to 
$V^{\otimes n}$, therefore $\DH$ can not be defined in $V$. 

Note also that any monomial of elements of $B$ in $v\in V$ has finite 
length. Since the deconcatenation of $v$ vanishes if repeated more then 
the length of $v$, the sum over $n$ in formula (\ref{DH-def}) is
indeed finite. 
For example, let us compute the simplest coproducts, for $x$, $y$ and $z$
in $B$:
\begin{eqnarray*}
\DH (x) &=& \sum (x\i1) \otimes (x\i2),\\
\DH (x,y) &=& \sum (x\i1,y\i1) \otimes (x\i2)(y\i2)
  + \sum (x\i1\cdot y\i1) \otimes (x\i2,y\i2),\\
\DH (x,y,z) &=& \sum (x\i1,y\i1,z\i1) \otimes
     (x\i2)(y\i2)(z\i2)
  + \sum (x\i1,y\i1\cdot z\i1) \otimes (x\i2)(y\i2,z\i2)
\\&&
  + \sum (x\i1\cdot y\i1,z\i1) \otimes (x\i2,y\i2)(z\i2)
  + \sum (x\i1\cdot y\i1 \cdot z\i1) \otimes 
       (x\i2,y\i2,z\i2).
\end{eqnarray*}
It is important to remark that, if $v=(x_{n_1},\dots,x_{n_l})$, 
then $\DH v$ is a sum of elements of $B^{\otimes k}$ tensorized by 
elements of $\T(B)^{\otimes k}$, for $k=1$ to $l$.
\medskip 

We now consider the dual Hopf algebra $\H^*=T(V^*)$, where $V^*$ is
the graded dual space $\T(B^*)$. 
Here $B^*$ is the dual bialgebra of $B$, with basis elements $x_m^*$ 
such that $\langle x_m^*\ |\ x_n \rangle=\delta_{m,n}$.
The product and coproduct in $B^*$ are defined in terms of those 
of $B$ by the standard duality 
\begin{align*}
\langle x^*\cdot y^*,z \rangle &= \langle x^* \otimes y^*, \delta
z\rangle, \\
\langle \delta x^*, y \otimes z \rangle &= \langle x^*, y\cdot
z\rangle.
\end{align*}
In $V^*=\T(B^*)$ we have $(x_1,\dots,x_n)^*=(x_1^*,\dots,x_n^*)$ and 
therefore 
\begin{align*}
\big\langle (x_1^*,\dots,x_n^*)\ |\ (y_1,\dots,y_m)\big\rangle 
&= \delta_{n,m} \langle x_1^*\ |\ y_1\rangle \cdots 
\langle x_n^*\ |\ y_n\rangle.
\end{align*}
As we said, $\H^*$ is a right-sided cofree-coassociative CHA, 
with deconcatenation coproduct $\D^*$ and product $\star$. 
In this section we describe the brace product induced by $\star$ 
on the set $V^*$ of primitive elements of $\H^*$. 
\medskip 

\begin{theorem}
The brace product $M_{1q}$ on $V^*=\T(B^*)$ which induces the 
CHA structure on $\H^*$ is given on the elements $v^*=(x^1,...,x^n)^*$ and 
$u_1^*,...,u_q^*$, for any $q>0$, by 
\begin{align*}
M_{1q}(v^*;u_1^*\cdots u_q^*) &= 0, \quad\mbox{if $n\neq q$}, 
\end{align*}
and, for $n=q$, by 
\begin{align*}
M_{1q}(v^*;u_1^*\cdots u_q^*) 
&= \big(x^1\i1 \cdot y_1^1\ ,...,\ x^1\i{k_1} \cdot y_1^{k_1}\ ,
...,\ x^q\i1 \cdot y_q^1\ ,...,\ x^q\i{k_q} \cdot y_q^{k_q}\big)^*, 
\end{align*}
where $u_i^*=(y_i^1,...,y_i^{k_i})^*$, for $i=1,...,q$. 
\end{theorem}

\begin{proof}
If $n\neq q$, we already remarked in the expression (\ref{DH-def}) 
of the coproduct that no term can occur. 

If $n=q$, consider $v=(x^1,...,x^q)\in B^{\otimes q}$ and 
$u_1,\dots,u_q$ in $B^{\otimes k_1}, \dots, B^{\otimes k_q}$,
respectively. 
Set $a=u_1\dots u_q \in V^{\otimes q}$. 
Then
\begin{eqnarray*}
M_{1q}(v^*;a^*) &=& \sum_k \sum_{w\in B^{\otimes k}}
\langle v^*\otimes a^*\ |\ \DH w \rangle w^* 
= \sum_k \sum_{w\in B^{\otimes k}} \sum 
\langle v^*\otimes a^*\ |\ w\i1 \otimes w\i2 \rangle w^*,
\end{eqnarray*}
where the sums run over the generators $w$ of $B^{\otimes k}$, 
that is, elements of the form $w=(z^1,\dots,z^k)$ where $z^i$ are
generators of $B$, and where $\sum w\i1 \otimes w\i2 = \d(w)$ 
denotes the coproduct of $B$. 

As a consequence, the term corresponding to $w$ appears in 
$M_{1q}(v^*;a^*)$ if $w\i1=v$ and $w\i2=a$. 
The equality $w\i1=v=(x^1,\dots,x^q)$ can be rewritten
\begin{eqnarray*}
z\i1^1\cdot{\dots}\cdot z\i1^{k_1} &=& x^1,\\
z\i1^{k_1+1}\cdot{\dots}\cdot z\i1^{k_1+k_2} &=& x^2,\\
&\dots &\\
z\i1^{k_1+\dots+k_{q-1}+1}\cdot{\dots}\cdot z\i1^{k} &=& x^q. 
\end{eqnarray*}

The equality $w\i2=a=u_1\dots u_q$ can be rewritten
\begin{eqnarray*}
(z\i2^1,\dots,z\i2^{k_1}) &=& u_1,\\
(z\i2^{k_1+1},\dots, z\i2^{k_1+k_2}) &=& u_2,\\
&\dots &\\
(z\i2^{k_1+\dots+k_{q-1}+1},\dots, z\i2^{k}) &=& u_q.
\end{eqnarray*}
The second equality means that $k=k_1+\dots+k_q$. 
%\todo{Check this last sentence, I dind't understand: is it $q$ or $k$?}
Because of this special form, each block can be treated 
independently. 

We consider the first block and we write $k=k_1$ and
$u_1=(y_1,\dots,y_k)$. 
We have now, for the first block and after putting $x=x^1$, 
\begin{eqnarray*}
X &=& \sum_{z^1,\dots,z^k} 
  \langle (x^*)\otimes (y_1,\dots,y_k)^* ,
  (z\i1^1\cdot{\dots}\cdot z\i1^{k})\otimes
  (z\i2^1,\dots,z\i2^{k})\rangle (z^1,\dots,z^k)^*
\\&=&
 \sum_{z^1,\dots,z^k} 
  \langle (x^*), (z\i1^1\cdot{\dots}\cdot z\i1^{k})\rangle
  \langle (y_1,\dots,y_k)^* ,
  (z\i2^1,\dots,z\i2^{k})\rangle (z^1,\dots,z^k)^*
\\&=&
 \sum_{z^1,\dots,z^k} 
  \langle x\i1^*, z\i1^1\rangle
  \dots
  \langle x\i{k}^*, z\i1^{k}\rangle
  \langle y_1^*,z\i2^1\rangle \dots
  \langle y_k^* , z\i2^{k}\rangle (z^1,\dots,z^k)^*,
\end{eqnarray*}
where we used the fact that the coproduct in $B^*$ is the dual of the 
product in $B$. Now we use the fact that the product in $B^*$ is the 
dual of the coproduct in $B$, and obtain
\begin{eqnarray*}
X &=& 
 \sum_{z^1,\dots,z^k} 
  \langle x\i1^*\cdot y_1^*, z^1\rangle
  \dots
  \langle x\i{k}^*\cdot y_k^*, z^k\rangle
  (z^1,\dots,z^k)^*
\\&=&
\sum 
  (x\i1\cdot y_1, \dots ,x\i{k}\cdot y_k)^*.
\end{eqnarray*}
The same is done for each block and we finally find
\begin{eqnarray*}
M_{1q}(v^*;a^*) &=& 
(x\i1^1\cdot y^1_1, \dots ,x\i{k_1}^1\cdot y^1_{k_1},
\dots,
x\i1^q\cdot y^q_1, \dots ,x\i{k}^q\cdot y^q_{k_q})^*. 
\end{eqnarray*}
\end{proof}

\begin{examples}
Examples of $M_{11}$ are: 
\begin{align*}
M_{11}\big((x^*);(y^*)\big) &= (x^*\cdot y^*), \\ 
M_{11}\big((x^*);(y^*,z^*)\big) &= 
   \sum (x^*\i1 \cdot y^*, x^*\i2 \cdot z^*), \\ 
M_{11}\big((x^*); (y^*_1, \dots,y^*_n) \big) &= 
   (x^*\i1 \cdot y^*_1, \dots, x^*\i{n} \cdot y^*_n).
\end{align*}

Examples of $M_{12}$ are: 
\begin{align*}
M_{12}\big((x,y)^*;(s)^*(t)^*\big)&=(x\cdot s,y\cdot t)^*,\\
M_{12}\big((x,y)^*;(s)^*(t_1,t_2)^*\big)&=
   \sum (x\cdot s, y\i1\cdot t_1,y\i2\cdot t_2)^*,\\
M_{12}\big((x,y)^*;(s_1,s_2)^*(t)^*\big)&=
   \sum (x\i1\cdot s_1, x\i2\cdot s_2,y\cdot t)^*, \\ 
M_{12}\big((x^*,y^*);(s^*_1, \dots,s^*_m)(t^*_1, \dots,t^*_n)\big) &= 
   \sum (x^*\i1 \cdot s^*_1, \dots, x^*\i{m} \cdot s^*_m,
   y^*\i1 \cdot t^*_1, \dots, y^*\i{n} \cdot t^*_n).
\end{align*}
\end{examples}

%%%%%%%%%%%%%%%%%%%%%%%%%%%%%%%%%%%%%%%%%%%%%%%%%%%%%

\end{document}